\documentclass[11pt]{article}
\usepackage{amsfonts, amsbsy, amssymb, eucal}
\usepackage{graphicx}

\begin{document}


\def\L{{\mathbb L}}
\def\B{{\mathbb B}}
\def\R{{\mathbb R}}
\def\P{{\mathbb P}}
\def\Rn{{\mathbb R}^n}
\def\Rnp{{\mathbb R}^{n+1}}
\def\T{{\mathbb T}}
\def\N{{\mathbb N}}
\def\M{{\mathbb M}}
\def\C{{\mathbb C}}
\def\Z{{\mathbb Z}}
\def\Q{{\mathbb Q}}
\newcommand{\ld}{,\ldots,}
\newcommand{\spr}[2]{\langle#1,#2\rangle}
\newcommand{\avr}[1]{\langle#1\rangle}

\newtheorem{definition}{Definition}
\newtheorem{assumption}{Assumption}
\newtheorem{theorem}{Theorem}
\newtheorem{lemma}{Lemma}[theorem]
\newtheorem{corollary}{Corollary}[theorem]

\title{On inverse problem of dynamics}
\author{M. Rudnev\thanks{Contact address: Department of Mathematics,
University of Bristol, University Walk, Bristol BS8 1TW, UK;\newline
e-mail: {\tt m.rudnev@bris.ac.uk}} $\,$and V. Ten\thanks{Contact
address: Department of Mathematics, University of Bristol,
University Walk, Bristol BS8 1TW, UK;\newline  e-mail: {\tt
v.ten@bris.ac.uk} }}
\date{February 13, 2004}
\maketitle

\begin{abstract}
We study the question whether for a natural Hamiltonian system on a
two-dimensional compact configuration manifold, a single trajectory
of sufficiently high energy is almost surely enough to reconstruct a
real analytic potential.
\end{abstract}

\def\proclaim#1{{\bf #1.}\bgroup\it}
\def\endproclaim{\egroup}
\def\leqslant{\leq}
\def\Bbb#1{{\bf #1}}

\medskip
\medskip
\noindent {AMS subject classification 37J05, 70H12}

\medskip
\medskip
\noindent Consider a compact configuration manifold $M^n$
(essentially we deal with the dimension $n=2$) and a natural
Hamiltonian system thereon, with Hamiltonian
$$
H(p,q)=\langle p,p\rangle_q+U(q),\;\;\;(p,q)\in T^*M^n.
$$
Above, $\langle \cdot,\cdot\rangle_q$ is a Riemannian metric on
$M^n$ and $U$ a potential. The direct problem of dynamics on $M^n$
is finding the trajectory $q(t)\subset M^n,$ with initial conditions
$q(0)=q_0$ and $\dot{q}(0)=v_0$, moving in the known force field
$f(q)=-\nabla_q U(q)$ on $M^n$, where the gradient $\nabla_q$ has
been associated with the metric $\langle \cdot,\cdot\rangle_q$.

Let us call the inverse problem of dynamics the problem of
reconstruction of the potential by observing the system's
trajectories $q(t)$. The first problem of this type was explored in
Newton's Principia, in quest for a physical law determining the
planetary motion compatible with observational data\footnote{Note
that the term ``inverse problem of mechanics'' has also been used to
address the problem of deciding whether a given system of second
order ODEs on $M^n$ has a Lagrangian, see e.g. \cite{LR}.}. In the
general case, knowledge of infinitely many trajectories is required
to completely solve the problem. In this note we show that in the
special case when $M^n$ is two-dimensional, compact and
topologically non-trivial, a single trajectory with sufficiently
large energy would almost surely suffice to reconstruct the
potential.

In the sequel, we assume that $M^n$ as well as all the quantities
involved are real-analytic. Also suppose, there is an a-priori
estimate $|U(q)|<C_0, \,\forall q\in M^n$, and we consider only the
trajectories $q(t)$ with total energy $E\geq C_0$, so $M^n$
coincides with the domain of possible motions.

\medskip
\noindent \proclaim{Theorem} Let $n=2,$ suppose $M^2$ is not
diffeomorphic to $S^2$ or $\R P^2$. Almost every trajectory
$q(t),\,t\geq0$, with energy $E\geq C_0$, suffices to reconstruct
the potential $U$ as a real-analytic function on $M^2$.
\endproclaim

\medskip
\noindent Let us recall the definition of a {\em key set}, or set of
uniqueness, see e.g. \cite{DS}.

\proclaim{Definition} Let $D$ be a domain in $\R^n$ and
$C^\omega(D)$ the class of real-analytic functions in $D$. A set
$K\subset D$ is a key set if any $f\in C^\omega(D)$ vanishing
identically on $K$, vanishes identically on $D$.
\endproclaim

\medskip
\noindent The theorem will follow from the following lemma.

\proclaim{Lemma} If a real-analytic dynamical system
$A:\;\dot{x}=F(x)$ on a compact phase space $P$,  is non-singular
(i.e. for no $x\in P,\,F(x)=0$) and has a positive measure of closed
orbits, then all its orbits are closed.
\endproclaim

\subsubsection*{Proof of the lemma} Since the manifold $P$ is compact
and the vector field $F$ is non-singular, there exists a Riemannian
metric $g$ on $P$, such that in this metric the vector field $F$ has
unit length at every $x\in P$. Furthermore, by compactness of $P$,
the curvature (associated with $g$) of integral trajectories of $A$
is bounded from above, and therefore there exists some $T_m>0$, such
that any periodic orbit of $A$ has period not smaller than $T_m$.

Let us partition the range $[T_m,\infty)$ of possible periods
(henceforth periods stand for minimum periods) for closed orbits on
intervals of some small length $\delta_1$ to be specified. Let $I_k
= [T_m+k\delta_1,T_m+(k+1)\delta_1)\equiv[T_k,T_{k+1})$, for
$k=0,1,\ldots.$ Let $\Gamma_k$ be the set of all closed orbits,
whose periods lie in $I_k$. Then for some $k=k_*$ the set
$\Gamma_*=\Gamma_{k_*}$, considered as a subset of $P$, has positive
measure. (We refer to the sets $\Gamma$ either as point sets or sets
of orbits, depending on the context. As there are only
measure-theoretical considerations involved, this should not cause
confusion.)

Let $D_{\delta_2}(x)$ be a codimension one disk in $P$, centered at
some $x\in P$, with radius $\delta_2$ and perpendicular (in the
sense of metric $g$) to the vector field $F(x)$ at the point $x$.
Let $\delta_2$ be small enough, so that the vector field is
transversal to $D_{\delta_2}(x)$ at every point of the disk.
Clearly, $\delta_2$ can be taken as a universal constant independent
of $x$.

Let $x_0$ be a Lebesgue point of the set $\Gamma_*$. Recall that at
a Lebesgue point, the density of the set is one. Let $\gamma_0$ be
the closed trajectory passing through $x_0$, so
$\gamma_0\in\Gamma_*$. Take a disk $D_{\delta_2}(x_0).$

Then there is a well defined analytic Poincar\'e map $S$ from a disk
$D_{\epsilon}(x_0)$ contained in  $D_{\delta_2}(x_0)$, where
$\epsilon<{1\over C_1 T_{k_*}}\delta_2$, for some
$C_1=C_1(\delta_1)$, and $x_0$ is a fixed point of $S$. On the disk
$D_{\epsilon}(x_0),$ points which are initial conditions for orbits
from $\Gamma_*$ form a set of positive measure, as $x_0$ is a
Lebesgue point. Besides, the quantity $\delta_1 < T_m$ can be chosen
small enough to ensure that any $x\in D_{\epsilon}(x_0)\cap
\Gamma_*$ is also a fixed point of the map $S$.

The union of all $x\in D_{\epsilon}(x_0)\cap \Gamma_*$ is a positive
measure subset of $D_{\epsilon}(x_0)$, and hence is a key set (see
\cite{Ten} for the proof that every set of positive measure is a key
set). Therefore, every point of $D_{\epsilon}(x_0)$ is an
equilibrium of the map $S$, and hence an initial condition for a
periodic orbit of $A$. Let $\Gamma_\epsilon$ be the union of all
such orbits, with initial conditions in $D_{\epsilon}(x_0)$. Let
$\Gamma\subseteq P$ be the maximum connected open set, which
contains $\Gamma_\epsilon$ and is a union of periodic orbits.

To complete the proof, let us show that the set $\Gamma$ does not
have a boundary, i.e. $\Gamma=P$. To show it, we use the following
Gronwall type estimate.

\proclaim{Proposition} Let $\phi(s)\subset \Gamma$ be a curve of
length $L$, where $s$ is a natural parameter with respect to the
metric $g$. Then for all $s\in [0,L]$ and some absolute constant
$C_2$,
$$
T(\phi(s))\leq T(\phi(0)) e^{C_2 L},
$$
where $T(\phi(s))$ is the period of the closed orbit passing through
the point $\phi(s)$.
\endproclaim

Indeed, the proposition follows immediately from the following
infinitesimal estimate: for some $C_2$,
$$
\left|{d\over ds}T(\phi(s))\right|\leq C_2 T(\phi(s)).
$$

Returning to the proof of the lemma, suppose the boundary
$\partial\Gamma$ is non-empty. As $\partial\Gamma$ is a compact set,
the distance (in the sense of metric $g$) between $\partial\Gamma$
and the above mentioned point $x_0$ attains its minimum at some
point $y\in\partial\Gamma$. Connect $x_0$ and $y$ by a geodesic
segment. Let the latter segment have length $L$; clearly all its
points, except $y$ belong to $\Gamma$. Let $\gamma_y$ be the
trajectory of $A$ with initial condition $y$.

For any point $x_1\neq y$ on the above geodesic segment, there is a
uniform bound for the period of the corresponding closed orbit, by
the proposition. Hence, for any such $x_1$, there exists a uniform
$\varepsilon$ (one can take $\varepsilon = \epsilon e^{-C_2 L}$,
where $\epsilon$ has been defined earlier) such that an analytic
Poincar\'e map can be defined in exactly the same way as $S$ above,
but now with the domain $D_\varepsilon(x_1)$. Choosing $x_1$ such
that the intersection $\gamma_y\cap D_\varepsilon(x_1)$ is not empty
and repeating the key set argument leads to contradiction: all
orbits in some tubular neighborhood of $\gamma_y,$ including
$\gamma_y$ itself, are periodic. $Q.E.D.$

\subsubsection*{Proof of the theorem} Consider a
randomly chosen trajectory $\gamma$ on some energy level
$H^{-1}(E),\,E\geq C_0$, which is obviously a non-critical level.

According to the lemma, either (i) all the trajectories on
$H^{-1}(E)$ are periodic, or (ii) a randomly chosen initial
condition $(p_0,q_0)\in H^{-1}(E)$ results in a phase trajectory of
infinite length almost surely.

The former case (i) may occur only if $M^2$ is  a so-called
$P$-manifold. Indeed, according to the Maupertuis principle, the
phase trajectories of motions with total energy  $E$ project onto
$M^2$ as geodesic lines of the corresponding (Riemannian) Jacobi
metric. $P$-manifolds are Riemannian manifolds, all whose geodesics
are closed, see \cite{Be}. Topological properties of $P$-manifolds
are characterized in great detail in various dimensions, within the
framework of the Bott-Samelson theorem. In our (simplest possible)
case, it is easy to see that $M^2$ can only be diffeomorphic to
either $S^2$ or $\R P^2$.

Indeed, the proof of the lemma implies that all the (closed) phase
trajectories on the energy level $E$ are homotopic to one another.
Then their images on $M^2$ (under natural projection) are also
homotopic. On the other hand, if $M^2$ is different from $S^2$ or
$\R P^2$, the number of generators for its fundamental group equals
at least two. As for any Riemannian metric there are closed
geodesics in each free homotopy class, we would have a
contradiction, unless $M^2$ is $S^2$ or $\R P^2$.

In the case (ii), let $q(t)=(q_1(t),q_2(t))$ be a randomly chosen
trajectory: it almost surely has infinite length. Clearly, we can
easily derive the gradient of $U$  at every point of the trajectory
from the Euler-Lagrange equations.

Thus to complete the proof of the theorem, let us show that any
non-closed trajectory is a key set. Consider orbit segments
$\{q(t),\,t\in[k,k+1)\}$, $k=0,1,\ldots$ in $M^2$. (Note that time
can always be scaled to ensure that each segment is a simple curve
in $M^2$, or shorter time intervals can be considered.) As $M^2$ is
compact, there is a limit point $q_*$ of the point sequence
$\{q(k+1/2)\}$. Consider a sufficiently small circle centered at
$q_*$. There are two options. Either the circle intersects the
trajectory $q(t)$ at infinitely many points, or at some point on the
circle, the trajectory $q(t)$ intersects itself (transversely as it
is a geodesic) infinitely many times. In the later case, take the
point of infinite self-intersection for $q_*$, otherwise leave $q_*$
as it is. In either case, there exists a point $q_*$, with the
property that any sufficiently small circle centered at $q_*$ is
intersected by the trajectory $q(t)$ infinitely many times.

Therefore, the force $f$ and the potential $U$ can be uniquely
reconstructed on any sufficiently small circle centered at $q_*$ (an
infinite point set on a circle is a key set for the circle), and
therefore in some neighborhood of $q_*$, and hence on the whole
configuration manifold $M^2$. $Q.E.D.$

\medskip
\noindent In conclusion, let us make several remarks.
\begin{enumerate}
\item The theorem can restated in terms of analytic
geodesic flows on compact Riemannian 2-manifolds. Namely, if $M^2$
is not diffeomorphic to a sphere or real projective plane, a
randomly chosen geodesic suffices to reconstruct the metric, almost
surely. Indeed, any Riemannian metric is locally conformal to the
Euclidean one, i.e. can be locally associated with Hamiltonian
$H(p,q)=e^{\rho(q_1,q_2)}(p_1^2+p_2^2)$. Our theorem enables one to
reconstruct the real-analytic function $\rho(q_1,q_2)$ locally near
$q_*$ from the Hamilton equations, with subsequent analytic
continuation to get the metric globally on $M^2$.

\item The theorem is essentially two-dimensional, as for $n\geq3$ the
fact that a random trajectory has infinite length does not suffice
to reconstruct the potential. Consider for instance the Euler top,
where $M^n=SO(3)$. In this case, the phase space is foliated by
invariant two-tori, where the trajectories are in general
conditionally periodic. Clearly, a projection of a single invariant
two-torus onto the three-dimensional configuration space is not a
key set.

\item The lemma does not apply to the case of invariant tori of dimension
higher than one. Indeed, a particular case of KAM theorem tells one
that a sufficiently small perturbation of a non-degenerate
Liouville-integrable Hamiltonian system in $T^*M^n$ yields a
positive measure set of invariant $n$-tori, which do not fill the
whole energy surface however.

\item
In the special case $M^2=\R P^2$, the only possibility to have all
the geodesic closed is the standard metric, by the theorem of L.
Green (\cite{G}). The case $M^2=S^2$ has been a subject of extensive
research for over 100 years, arguably beginning with the doctoral
thesis of O. Zoll (\cite{Zo}). The reader is referred to the
excellent book by L. Besse (\cite{Be}), which gives the issue a
thorough treatment.

\item
The condition $E\geq C_0$ also seems unavoidable, as for small
energies the domain of possible motions can be a disk, with the
Jacobi metric degenerate on the boundary, in which case one may
expect the scenario, similar to the case $M^2=S^2$.

\item
Observe that in the exceptional case when $M^2$ is a $P$-manifold,
the geodesic flow thereon is completely integrable (see e.g.
\cite{BJ} for the proof of this fact). Hence our theorem implies
that if $n=2$, it is sufficient for restoration of the potential
almost surely from a single trajectory that the system possess no
other analytic integrals of motion but energy. It seems likely that
in the latter weaker formulation the theorem should extend to the
case $n>2$, however we do not know how to prove it.

\end{enumerate}
\bibliographystyle{unsrt}

\end{document}